\documentclass[11pt]{amsart}
\usepackage{amsmath}
\usepackage{amssymb}
\usepackage{amsthm}
\usepackage{graphicx}

\usepackage{color}

\newcommand{\Q}{\mathbb{Q}}
\newcommand{\Z}{\mathbb{Z}}
\newcommand{\F}{\mathbb{F}}

\newcommand{\Mat}{\mathrm{Mat}}
\newcommand{\GL}{\mathrm{GL}}
\newtheorem{lem}{Lemma}
\newtheorem{theorem}[lem]{Theorem}
\theoremstyle{definition}
\newtheorem{defi}[lem]{Definition}
\theoremstyle{remark}
\newtheorem{rem}[lem]{Remark}

\newcommand{\la}{\langle}
\newcommand{\ra}{\rangle}
\begin{document}
\title{Presentations of matrix rings}
\author{Martin Kassabov}
\thanks{The author was supported in part by the NSF grants
DMS~0600244 and DMS~0900932.}

\maketitle

Recently, there has been
a significant interest in the combinatorial properties the ring
of $n\times n$ matrices.
The aim of this %short
note is to describe a short (may be the shortest possible)
presentation of the matrix ring $\Mat_n(\Z)$. This presentation
%,described below,
is
significantly shorter than the previously known ones, see~\cite{SP}.
\iffalse
We hope that this presentation will have applications in computational
ring theory and might also be used to build crypto-systems based on
combinatorial problems in rings.
\fi

Surprisingly, the number of relations in the presentation does
not depend on the size of the matrices and all relations have
relatively simple form. In contrast, the similar statement for
the groups $\GL_n(\Z)$ is significantly more difficult to prove
and the presentations have more relations, see~\cite{GKKL1,GKKL3}.

\begin{theorem}
\label{th1}
The ring%
\footnote{All rings in this paper are associative and contain a unit element.
Also all presentations are in the category of unitial associative rings.}
of $n\times n$ matrices over $\Z$, for $n \geq 2$, has a presentation with
$2$ generators and 3 relations
$$
\Mat_n(\Z) = \la x,y \mid x^n = y^n=0, xy + y^{n-1}x^{n-1} = 1 \ra.
$$
\end{theorem}
\begin{proof}
Let $R$ denote the ring defined by the above presentation.
A homomorphism from the free associative ring to $\Mat_n(\Z)$ given by
$$
x \to X= \sum e_{i,i+1}
\quad
y \to Y= \sum e_{i,i-1}
$$
factors through the ring $R$.%
\footnote{As usual, $e_{i,j}$ denotes the elementary matrix with $1$ at $i,j$-th place
and zeroes everywhere else.}
The first two relations are satisfied
because both $X$ and $Y$ are nilpotent matrices and the third follows
form a direct computation.
\iffalse
shows that they satisfy the relation
$$
XY + Y^{n-1}X^{n-1} = 1.
$$
\fi
Moreover, this
homomorphism is surjective because
the matrices $X$ and $Y$ generate
$\Mat_n(\Z)$ as a ring.

Thus, it remains to prove that the ring $R$ is not too big. Our first step
is to find other relations, which are satisfied in the ring $R$.

\begin{lem}
\label{lem1}
For any non-negative integers $k,l,m \geq 0$, such that $m\geq l$, we have
$$
x^ky^lx^m = \left\{
\begin{array}{lr}
y^{l-k}x^k  & \mbox{if } l \geq k \\
x^{k+m-l}  & \mbox{if } l \leq k
\end{array} \right.
$$
and
$$
y^mx^ly^k = \left\{
\begin{array}{lr}
y^{m}x^{l-k}  & \mbox{if } l \geq k \\
y^{k+m-l}  & \mbox{if } l \leq k
\end{array} \right. .
$$
\end{lem}
\begin{proof}
The proof uses induction on $l$. The base case $l=0$ is trivial.
If $k=0$ again there is noting to prove. Thus, without loss of generality
we can assume that $k,l\geq 1$.
$$
x^{k}y^l x^m = x^{k-1} (xy) y^{l-1} x^{m}=
x^{k-1} (1 - y^{n-1}x^{n-1}) y^{l-1} x^{m} =
$$
$$
=x^{k-1}y^{l-1}x^m - x^{k-1}y^{n-1}\left( x^{n-1}y^{l-1}x^{m}\right).
$$
By the induction assumption the part of the second term in the brackets
is equal to $x^{n+m-l}=0$, since $m \geq l$. Thus, the second term
vanishes.
Another application of the induction assumption shows that
$$
x^{k-1}y^{l-1}x^{m} =\left\{
\begin{array}{lr}
y^{l-k}x^{m} & \mbox{if } l \geq k \\
x^{k+m-l}   & \mbox{if } l \leq k
\end{array} \right.
$$
and completes the proof of the first part of the lemma.

The proof of the second part uses similar induction.
Alternatively, one can use that the transformation
$x\to y$ and $y\to x$ extends to an %involution
anti-automorphism $\sigma$ of $R$ and
apply $\sigma$ to the identities from the first part.
%notice
%that the identities in second part are just images of the ones in the first part under this
%anti-isomorphism.
\end{proof}

\begin{lem}
\label{lem2}
The ring $R$ is generated, as an additive group, by the elements
$y^ix^j$ for $0 \leq i,j <n$.
\end{lem}
\begin{proof}
Let $T$ denote the additive subgroup of $R$ generated by the elements $y^ix^j$.
% Since $1\in T$ in order to prove that $T=R$,
It suffices to show that $T$ is closed under, both left and right, multiplication
by $x$ and $y$, because $1=x^0y^0 \in T$.
The two relations $x^n=y^n=0$ imply that
$T$ is closed under left multiplication by $y$
and right multiplication by $x$. Thus, it remains to show that $xT, Ty \subseteq T$.
The element $x. y^ix^j$ is clearly in $T$ if $i=0$. If $i \geq 1$ we have
$$
xy^ix^j =(xy) y^{i-1} x^j =
\left( 1 - y^{n-1}x^{n-1}\right) y^{i-1} x^j =
$$
$$
= y^{i-1} x^j -\left(y^{n-1}x^{n-1}  y^{i-1} \right)x^j =
y^{i-1} x^j - y^{n-1} x^{n-i} x^j \in T,
$$
where the last equality uses Lemma~\ref{lem1}.

Similarly, $y^ix^j.y$ is in $T$ if $j=0$ and if $j\geq 1$ we have
$$
y^ix^jy = y^{i} x^{j-1}(xy)=
y^{i} x^{j-1}\left( 1 - y^{n-1}x^{n-1}\right) =
$$
$$
= y^{i} x^{j-1} - y^{i} \left(x^{j-1}y^{n-1}x^{n-1} \right) =
y^{i-1} x^j - y^{i} y^{n-j} x^{n-1} \in T.
$$
\end{proof}

Lemma~\ref{lem2} together with the surjection for $R$ to
$\Mat_n(\Z)$ is sufficient to show to $R$ is isomorphic to
the matrix ring, but one can build the isomorphism directly:

\begin{defi}
Let $a_{i,j}$, for $0\leq i,j < n$, denote the elements
$$
a_{i,j} = y^ix^j - y^{i+1}x^{j+1}.
$$
\end{defi}

\begin{lem}
\label{lem3}
We have that
$$
a_{i,j}x = a_{i,j+1}
\quad
a_{i,j}y = a_{i,j-1}
\quad
xa_{i,j} = a_{i-1,j}
\quad
ya_{i,j} = a_{i+1,j}.
$$
Here, we assume that $a_{i,j} = 0$ if either $i$ or $j$ is outside the interval $[0,n-1]$.
\end{lem}
\begin{proof}
Two of the identities follow directly from the definition of the elements $a_{i,j}$.
The only non-trivial ones are $a_{i,j}y = a_{i,j-1}$ and $xa_{i,j} = a_{i-1,j}$.
By the definition of the element $a_{i,j}$ we have% Let us prove the first one:
$$
a_{i,j}y = y^ix^jy - y^{i+1}x^{j+1} y
$$
If $j=0$ the right side is equal to $y^{i+1}  - y^{i+1}xy = y^{i+1}y^{n-1}x^{n-1} = 0$.
Otherwise, we can use the proof of Lemma~\ref{lem2}:
$$
=
\left(y^i x^{j-1} - y^{n+i-j}x^{n-1} \right) -
\left(y^{i+1} x^{j} - y^{n+i-j}x^{n-1} \right)=
$$
$$
=y^i x^{j-1} - y^{i+1} x^{j} = a_{i,j-1}.
$$
The proof of the second relation $xa_{i,j} = a_{i-1,j}$ is similar.
\end{proof}

\begin{lem}
\label{lem4}
The product $a_{i,j} a_{p,q}$ is equal to $0$ if $j\not = p$ and is equal to
$a_{i,q}$ if $j=p$.
Moreover we have
$$
1=\sum a_{i,i}
\quad
x=\sum a_{i,i+1}
\quad
y=\sum a_{i+1,i}.
$$
%Thus, the map $a_{i,j} \to e_{i,j}$ is an isomorphism between the ring $R$ and $\Mat_n(\Z)$.
\end{lem}
\begin{proof}
These equalities follows directly from
Lemma~\ref{lem3} and the definition of the elements $a_{ij}$.
\end{proof}

Thus, the map $a_{i,j} \to e_{i,j}$ extends to an isomorphism between the ring $R$ and $\Mat_n(\Z)$,
%By Lemmas~\ref{lem2} and~\ref{lem3} it follows that the elements $a_{i,j}$
%span the additive group of $R$ and they multiply as the elementary matrices
%in $\Mat_n(\Z)$,
which completes the proof of Theorem~\ref{th1}.
\end{proof}
\begin{rem}
From the isomorphism between $R$ and $\Mat_n(\Z)$ it follows that the elements
$x$ and $y$ also satisfy the relation
$$
yx + x^{n-1}y^{n-1}=1,
$$
i.e., the map $x\to y$ and $y\to x$ can be extended to an automorphism of $R$.
\end{rem}

\bigskip

The following variation of Theorem~\ref{th1} gives presentation of the
ring of matrices over $\Z/N\Z$:
\begin{theorem}
\label{th2}
For any integer $N$ the matrix ring
$\Mat_n(\Z/N\Z)$ has presentation
$$
\la x ,y \mid  x^n = y^n =0, xy + (N+1)y^{n-1}x^{n-1}=1 \ra.
$$
\end{theorem}
\begin{proof}
The argument is a slight modification of the proof of Theorem~\ref{th1}.
The map $x \to X$ and $y \to Y$ is a homomorphism onto
$\Mat_n(\Z/N\Z)$. The proof of injectivity follows the same outline ---
Lemmas~\ref{lem1}, \ref{lem2}, \ref{lem3}
and~\ref{lem4} still hold.
Finally one observes that $a_{n-1,n-1} = y^{n-1}x^{n-1}$
and $xy = \sum_{i=0}^{n-2} a_{i,i}$. Thus the last relation
in the presentation is equivalent to
$$
\sum_{i=0}^{n-1} a_{i,i} = 1 = xy + (N+1)y^{n-1}x^{n-1} =
$$
$$
=\sum_{i=0}^{n-2} a_{i,i}  + (N+1) a_{n-1,n-1}
$$
therefore $Na_{n-1,n-1}=0$, which implies that
$Na_{i,j}=0$ for all $i,j$, i.e., the additive group of $R$
has exponent $N$.
\end{proof}
\begin{rem}
In the case $N=p$ is a prime number, one can use that the matrix
algebra $\Mat_n(\F_p)$ is a cyclic algebra, thus it is possible to
obtain a presentation (as an algebra over $\F_p$) with $2$ generators
and $3$ relations. This leads to a presentation of the ring
$\Mat_n(\F_p)$ with $2$ generators and $4$~relations.
The presentation obtained using this approach
uses a presentation of the finite field $\F_{p^n}$, which involves
an irreducible polynomial of degree $n$ over $\F_p$. Thus,
the relations in this presentation will be more ``complicated'' than
the ones %in the presentations in
Theorem~\ref{th2}.
In some cases, one can modify the presentation of the
cyclic algebra and save one relation:%, see~\cite{Gu-private}:%
\footnote{This presentation was found by Robert Guralnick~\cite{Gu-private}.}
$$
\Mat_p(\F_p) = \la x, y \mid y^p=1, x^p=x, xy = y(x+1) \ra.
$$
\end{rem}

\begin{rem}
In some sense the presentation in Theorem~\ref{th1} is a
variant of the presentation of cyclic algebra, where one uses
the nilpotent ring $\Z[x]/(x^n)$ instead of the maximal subfield,
however it is not completely clear what is the analog of the
``field'' automorphism in this picture.
\end{rem}

\medskip

One would like to say that the presentations in Theorems~\ref{th1}
and~\ref{th2} are the
simplest possible. Unfortunately, we were not able to prove that
presentations of the matrix rings with $2$ generators and
only $2$ relations do not exist.
The following result shows that there
isn't any
%does not exist a
presentation
of the matrix ring $\Mat_n(\Z)$ with a single relation:
\begin{theorem}
\label{th3}
The number of relations in any presentation of the matrix ring $\Mat_n(\Z)$
is at least equal to the number of generators.
\end{theorem}
\begin{proof}
The main idea of the proof is to ``translate'' the notion of the ``relation module''
from groups rings and use it to obtain a lower bound for the number of relations
in a presentation of ring, see~\cite{I}.

Let $ 0 \to I \to \Z\la S \ra \to R \to 0$ be a presentation of the ring $R$.
The quotient $I/I^2$ is called {\emph{relation module associated to this presentation}}
and is naturally a $R$ bi-module. It is clear that the projection of any generating set of the ideal $I$
to $I/I^2$ is a generating set of this %module as an
bi-module. Thus, the minimal
number of generators of relation module gives a lower bound for the number of generators of
the ideal $I$ and the number of relations in
a presentation of the ring $R$.

%In the case of groups the relation module corresponds to the maximal abelian extension.
%In the case of rings the role of abelian extension is played by radical extensions.

One way to construct a big quotient of the relation module
%\footnote{This construction often yields the torsion free
%part of the relation module.}
is the following:
Let $d$ be the size of the
generating set $\bar S$ and let $M$ be the free $R$ bi-module on $d$ generators $m_i$, i.e.,
$M \simeq (R \otimes R)^{\oplus d}$.
We can define a ring structure on the abelian group $R \oplus M$, where the multiplication
between elements of $R$ and $M$ is defined using left and right actions of $R$ on $M$ and
the product of any two elements in $M$ is equal to zero.
%Equivalently, one can consider the ring of matrices of the form
%$\left( \begin{array}{cc} R & M \\ )

For any generating set $\bar S$ of $R$ with $d$ elements one can define a subring $\tilde R$ of $R\oplus M$
generated by ``extensions'' of the generators in $\bar S$ by generators of the module $M$,
i.e., $\tilde s_i = \bar s_i +m_i$. It is easy to see %can be shown
that the relation module corresponding to this presentation maps surjectively onto
the intersection of $\tilde R$ with $M$.

In the special case of the matrix ring $\Mat_n(\Z)$, we have
the module isomorphism $M\simeq \Mat_n(L)$, where
$L$ is free abelian group on $n^2d$ generators (this follows form the isomorphism
of the bi-modules $\Mat_n(\Z)\otimes \Mat_n(\Z) \simeq \Mat_n(\Z^{\oplus n^2})$).
%The submodules of $M$ correspond to subgroups of the group $L$.
A long computation shows that the intersection of $\tilde R$ with $M$ is
isomorphic to $\Mat_n(\tilde L)$, where $\tilde L \subset L$ is a subgroup of rank at least
$n^2(d-1) + n > n^2(d-1)$ --- this bound does not depend on the images of the
generators $S$ in $\Mat_n(\Z)$.
This completes the proof of the Theorem, since the relation module can not be
generated by than $d-1$ elements.
\end{proof}

\begin{rem}
A more carefully computation shows that the relation module, corresponding to
the map $\Z\la x,y \ra \to \Mat_n(\Z)$ given by $x\to X$ and $y\to Y$, can
be generated by only two elements as an %$\Mat_n(\Z)$
bi-module (for example the elements $xy+y^{n-1}x^{n-1}-1$ and $xy^n + yx^n$
generate the relation module).
This suggests that it might be possible to ``combine'' the
two  relations $x^n=y^n=0$ into a single one and obtain a presentation of
$\Mat_n(\Z)$ with $2$ generators and only $2$ relations.
Unfortunately, the usual trick of combining such relations --- replacing them
with $x^n=y^n$ does not work, because the presentation
$$
\la x,y \mid x^n=y^n, xy + y^{n-1}x^{n-1} = 1 \ra
$$
defines a ring which surjects onto $\Mat_n\left(\Z[t]/(t^2)\right)$.
\end{rem}

\begin{rem}
One can view the proof of Theprem~\ref{th3} as an analog of
Gasch{\"u}tz' result~\cite{G,Gru},
which says that
the tensor product of the relation module for a finite group $G$ with $\Q$ is
isomorphic as $G$-module to $\Q[G]^{\oplus d-1} \oplus \Q$,
therefore
the relation module can not be generated by
less than $d$ elements.
\end{rem}

{\bf Acknowledgment:}
The presentations in this paper arise as a side-result of a project
about presentations of finite groups. The author wish to thank
his collaborators in that project Robert Guralnick, William Kantor and
Alex Lubotzky  for useful discussions.
%The author also grateful to IAS for the hospitality during the work on that project.

%\bibliographystyle{plain}
%\bibliography{matrix}

\begin{thebibliography}{1}

\bibitem{G}
Wolfgang Gasch{\"u}tz.
\newblock \"{U}ber modulare {D}arstellungen endlicher {G}ruppen, die von freien
  {G}ruppen induziert werden.
\newblock {\em Math. Z.}, 60:274--286, 1954.

\bibitem{Gru}
Karl~W. Gruenberg.
\newblock {\em Relation modules of finite groups}.
\newblock American Mathematical Society, Providence, R.I., 1976.
\newblock Conference Board of the Mathematical Sciences Regional Conference
  Series in Mathematics, No. 25.

\bibitem{Gu-private}
R.~M. Guralnick.
\newblock private communacation.

\bibitem{GKKL3}
R.~M. Guralnick, W.~M. Kantor, M.~Kassabov, and A.~Lubotzky.
\newblock Presentations of finite simple groups: a computational approach.
\newblock {\em J. European Math. Soc., to appear}.

\bibitem{GKKL1}
R.~M. Guralnick, W.~M. Kantor, M.~Kassabov, and A.~Lubotzky.
\newblock Presentations of finite simple groups: a quantitative approach.
\newblock {\em J. Amer. Math. Soc.}, 21(3):711--774, 2008.

\bibitem{I}
S.~V. Ivanov.
\newblock Relation modules and relation bimodules of groups, semigroups and
  associative algebras.
\newblock {\em Internat. J. Algebra Comput.}, 1(1):89--114, 1991.

\bibitem{SP}
B.~V. Petrenko and S.~N. Sidki.
\newblock On pairs of matrices generating matrix rings and their presentations.
\newblock {\em J. Algebra}, 310(1):15--40, 2007.

\end{thebibliography}

%\bigskip

Martin Kassabov,
Department of Mathematics,
Cornell University,
Malott Hall,
Ithaca, NY 14853-4201, USA

email: {\tt kassabov@math.cornell.edu}

\iffalse
\begin{minipage}[t]{3.5 in}
Martin Kassabov\\
Department of Mathematics, \\
Cornell University, \\
Malott Hall,\\
Ithaca, NY 14853-4201, USA\\
email: {\tt kassabov@math.cornell.edu}
\end{minipage}
\fi
\end{document}